\documentclass[a4paper,twoside]{article}
\usepackage{geometry}
\usepackage{fancyhdr}
\usepackage{amsmath}
\usepackage{amssymb}
\usepackage{amsthm}
\usepackage{enumerate}
\usepackage{mathrsfs}
\usepackage{xcolor}
\usepackage{indentfirst}

\newtheorem{thm}{Theorem}[section]
\newtheorem{lem}[thm]{Lemma}

\newtheorem{prop}[thm]{Proposition}
\newtheorem{rem}[thm]{Remark}
\newtheorem{defi}[thm]{Definition}
\numberwithin{equation}{section}
\begin{document}
	\begin{center}
		{\Large \bf Random attractors for the stochastic Nernst-Planck-Navier-Stokes system with multiplicative white noise}\\
		\vspace{0.5cm} {Yang-yang WU$^*$ }\\
		{\small
			Department of Mathematics,\\
			Nanjing University of Aeronautics and Astronautics,\\
			Nanjing, 211106, P.R. China
		}\\
		\vspace{0.5cm}
		{Gao-cheng YUE}\\
		{\small
			Department of Mathematics,\\
			Nanjing University of Aeronautics and Astronautics,\\
			Nanjing, 211106, P.R. China
		}
	\end{center}
	\footnote[0]{
		\hspace*{0.6cm}$^*$ Corresponding author.\\
		\hspace*{0.6cm}\textit{E-mail address}: wyynjtech@163.com}
	
	\begin{abstract}
		In this paper, we consider the 2D periodic stochastic Nernst-Planck-Navier-Stokes equations with body forces perturbed by multiplicative white noise. We first transform the stochastic Nernst-Planck-Navier-Stokes system into the deterministic system and address the problem of global well-posedness of the solution. Then, we generate a corresponding random dynamical system and dedicate to proving the existence of a compact random attractor for such random dynamical system. Furthermore, upper semicontinuity of the random attractor is established when the noise intensity approaches to zero. 
	\end{abstract}
	\indent
	\hspace{1.0cm}{\bf Keywords:} Electrodiffusion; Random attractor; Stochastic Nernst-Planck-Navier-Stokes equations; White noise.\\
	\indent
	\hspace{0.4cm}{\bf 2020 Mathematics Subject Classification:} 35B41; 35Q35; 35R60; 76D03.

	\section{\large Introduction}
	In this paper, we consider the Nernst-Planck-Navier-Stokes (NPNS) system with fluid body forces. We study two ionic species with equal diffusivities $D>0$ and valences $z_1=1$ and $z_2=-1$, respectively, which is the simplest setting for the ionic concentrations. Mathematically, the evolution of ionic concentrations $c_1$ and $c_2$ obeys the Nernst-Planck equations
	\begin{equation}\label{y1}
		\partial_{t}c_1+u\cdot\nabla c_1=D\Delta c_1 +D\nabla\cdot(c_1\nabla \Phi),
	\end{equation}
	\begin{equation}\label{y2}
		\partial_{t}c_2+u\cdot\nabla c_2=D\Delta c_2 -D\nabla\cdot(c_2\nabla \Phi).
	\end{equation} 
	The fluid velocity $u$ is described by the Navier-Stokes system and the divergence-free condition
	\begin{equation}\label{y3}
		\frac{du}{dt} + u\cdot\nabla u + \nabla p- \nu\Delta u  =-\rho \nabla \Phi +f,~\nabla \cdot u=0.
	\end{equation}
	The potential $\Phi$ is given by the Poisson equation
	\begin{equation}\label{y4}
		-\varepsilon_0 \Delta\Phi=\rho=c_1-c_2.
	\end{equation}
	Above $\rho$ denotes the charge density, $\nu>0$  is the kinematic viscosity of the fluid, $\varepsilon_0>0$ is a positive constant proportional to the square of the Debye length, $p$ denotes the pressure, $\Phi$ represents the electrical potential, and $f$ is a time independent, smooth and divergence-free body force.
	
	From \eqref{y1} to \eqref{y4}, the NPNS system is obtained as follows
	\begin{equation}\label{a1}
		\left\{
		\begin{aligned} 
			\displaystyle{\frac{du}{dt} + u\cdot\nabla u + \nabla p- \nu\Delta u  =-\rho \nabla \Phi +f},  \\
			\nabla \cdot u = 0,  \\
			\rho=c_1-c_2, \\
			-\varepsilon_0 \Delta\Phi=\rho, \\
			\partial_{t}c_1+u\cdot\nabla c_1=D\Delta c_1 +D\nabla\cdot(c_1\nabla \Phi), \\
			\partial_{t}c_2+u\cdot\nabla c_2=D\Delta c_2 -D\nabla\cdot(c_2\nabla \Phi). \\
		\end{aligned}
		\right.
	\end{equation}
	The unknowns are the velocity $u$ and the ionic concentrations $c_1$ and $c_2$.
	
	From mathematical point of perspective, the system we are discussing encompasses classical Nernst-Planck equations and Navier-Stokes equations. There has been a tremendous literature on the well-posedness of the Nernst-Planck \cite{E.F.W2023, M.J2021, M.J2022} and Navier-Stokes systems \cite{P.C1988}. Electrodiffusion in fluids, which is described by the NPNS equations \cite{I.R1990}, has a	variety of outstanding applications in neuroscience \cite{N.T1989}, semiconductors \cite{H.K1986}, water purification, ion separations \cite{A.J2016, H.J2018}, and ion selective membranes \cite{G.J1989}. Due to its importance in real world, the deterministic Nernst-Planck-Navier-Stokes equations have been studied extensively in different dimensions and situations. For example, it has been proved that the NPNS system has global and stable solutions with blocking and selective boundaries in 2D \cite{P.M2019}. In the case of mixed blocking and selective boundary conditions for the ionic concentrations and Robin boundary conditions for the electric potential, the global existence of three-dimensional strong solutions was shown in \cite{Lee.2022}. When both the ionic concentrations and the electrical potential are under the Dirichlet boundary conditions, the 3D strong solutions of the NPNS system are global if the fluid velocity is regular in \cite{P.M.F2021}. In 2D periodic domains, the authors explored long time dynamical behavior of solutions and showed the convergence of the solutions with the presence of body charges and body forces in \cite{E.M2021}. In the periodic boundary conditions, the authors considered the long-time dynamics of the NPNS system and established the existence of a global unique smooth solution, see \cite{P.M2024}. In addition, the existence of a global analytic solution for the forced NPNS system was obtained in \cite{E.M2022}. More Recently, the authors studied the global well-posedness and long-time dynamics of Nernst-Planck Darcy system and further proved the existence of a finite-dimensional global attractor in \cite{E.D2024}.
	
	However, in many practical circumstances, random effects play a key role that can not be ignored. These effects include turbulence, random fluctuations and uncertainties in fluid systems, which are essential for achieving an accurate estimate of dynamical systems in realistic world. There exist some literature about stochastic models and their properties that are related to our further study. In \cite{B.W2014}, the author studied the existence and uniqueness of random attractors for stochastic Reaction-Diffusion equations with multiplicative noise on unbounded domains and proved the upper semicontinuity of attractors. Moreover, Some results and models concerning the stochastic Navier-Stokes equations have been obtained in \cite{L.F2020, M.T2024, M.P2021, R.B2004}. 
	
	As is known to us, multiplicative noise provides an appropriate description for stochastic disturbance and has thus attracted extensive attention. There exist numerous results with respect to the stochastic PDEs driven by linear multiplicative noise in \cite{Z.B.M1998,Z.B.M2016,Z.D.B.G2023,BX.W2014,D.Z2022}. In the two-dimensional case, the stochastic Navier-Stokes equation driven by linear multiplicative white noise has been studied in \cite{S.T2010}, with solutions of the 2D Navier-Stokes equation generating a locally compact cocycle. More recently, the authors established the existence of a global pullback attractor and proved the upper semicontinuity for damped Navier-Stokes equations of 2D as well as 3D in \cite{K.M2023}. 
	
	To the best of our knowledge, there is still a scarcity of research about random attractors for stochastic NPNS system. Motivated by the literature above, we consider the stochastic Nernst-Planck-Navier-Stokes (SNPNS) equations perturbed by multiplicative white noise as follows:
	\begin{equation}
		\left\{
		\begin{aligned}\label{a2} 
			\displaystyle{\frac{du}{dt} + u\cdot\nabla u + \nabla p- \nu\Delta u  =-\rho \nabla \Phi +f+ \varepsilon u\circ\ \frac{dW(t)}{dt} },  \\
			\nabla \cdot u = 0,  \\
			\rho=c_1-c_2, \\
			-\varepsilon_0 \Delta\Phi=\rho, \\
			\partial_{t}c_1+u\cdot\nabla c_1=D\Delta c_1 +D\nabla\cdot(c_1\nabla \Phi), \\
			\partial_{t}c_2+u\cdot\nabla c_2=D\Delta c_2 -D\nabla\cdot(c_2\nabla \Phi). \\
		\end{aligned}
		\right.
	\end{equation}
	We address the SNPNS system in the two-dimensional periodic domain $\mathbb{T}^{2}$=$[0,2\pi]$$\times[0,2\pi]$ with periodic boundary conditions. Here, $\varepsilon > 0$ represents the noise intensity, $\circ$ means that the stochastic integral is understood in the sense of Stratonovich and $W(\cdot)$ is the two-sided Wiener process on the probability space   $(\Omega,\mathscr{F},P)$, where $\Omega =\{\omega \in C(\mathbb{R},\mathbb{R}) :\omega(0)=0\} $, $\mathscr{F}$ is the Borel sigma-algebra induced by the compact-open topology of $\Omega$, and $P$ is a Wiener measure.
	
	In this paper, our main goal is to prove the existence of a compact random attractor and establish upper semicontinuity of the random attractor. The upper semicontinuity of random attractors was introduced in \cite{T.J.J1998}, that is, random attractors $\mathscr{A_\varepsilon}(\omega)$ converge towards a deterministic global attractor $\mathscr{A}(\omega)$ when the noise intensity $\varepsilon$ approaches to zero. 
	
	The rest of this article is organized as follows. In section \ref{sec1}, we review some functional spaces and inequalities needed throughout the paper. Section \ref{sec2} contains basic concepts about random dynamical system and random attractors. Section \ref{sec3} is concerned with the existence and uniqueness of the solution, as well as the continuous dependence of the solution with respect to the initial data. We finally obtain a compact absorbing set and show the existence of a compact random attractor in $H$ in section \ref{sec4}. Section \ref{sec5} is devoted to the
	upper semicontinuity of such random attractor.

	\section{\large Preliminaries}\label{sec1}
	In this section, we will introduce some functional spaces, classical notations, and common inequalities that will be used in the further analysis. For more details of these topics, we refer to \cite{P.C1988}. 
	
	We denote by $L^p$ the Lebesgue spaces with the norm
	\begin{align*}
		\|u\|_{L^p}=(\int _{\mathbb{T}^2}|u(x)|^p dx)^{\frac{1}{p}},~ ~\mathrm{if}~ p \in [1,\infty)
	\end{align*}
	and
	\begin{align*}
		\|u\|_{L^\infty}=\mathop {sup }\limits_{x \in \mathbb{T}^2}|u(x)|,~ ~\mathrm{if}~ p =\infty.
	\end{align*}
	We also define the following general abstract space
	$$\displaystyle{\tilde{V} = \left \{u \in \left (C_{0}^{\infty}({\mathbb{T}^{2}})\right )^{2}\ : \ \mathrm{div}\,u = 0\right \},}$$
	where $C_{0}^{\infty}({\mathbb{T}^{2}})$ denotes the space of all infinitely differentiable functions with compact support in $\mathbb{T}^{2}$. Let $H$ denote the closure of $\tilde{V}$ in $(L^2(\mathbb{T}^{2}))^2$ with the norm $\left \vert \cdot \right \vert$, and the inner product $(\cdot , \cdot)$, where for $u,\upsilon\in(L^2(\mathbb{T}^{2}))^2$, we set
	$$\displaystyle{(u,\upsilon) =\sum_{i=1}^{2}\int_{\mathbb{T}^{2}}u_{i}(x)\upsilon_{i}(x)\,dx.}$$
	Moreover, let $V$ denote the closure of $\tilde{V}$ in $(H_0^1(\mathbb{T}^{2}))^2$ with the norm $\|\cdot \|$, and the associated scalar product $((\cdot , \cdot))$, where for $u,\upsilon\in(H_0^1(\mathbb{T}^{2}))^2$, we have
	$$\displaystyle{((u,\upsilon)) =\sum_{ i,j=1}^{2}\int_{\mathbb{T}^{2}}\frac{\partial u_{j}} {\partial x_{i}} \frac{\partial \upsilon_{j}} {\partial x_{i}}\,dx.}$$
	It follows that $V \subset H \equiv H^{{\prime}}\subset V ^{{\prime}}$, where the injections are continuous and dense.
	Let $\left \langle \cdot,\cdot \right \rangle$ denote the duality pairing between $V$ and $V ^{{\prime}}$. 
	
	We recall the Stokes operator $A:V\rightarrow V^{\prime}$ by
	\begin{align*}
	 \langle Au,\upsilon\rangle =\int_{\mathbb{T}^{2}}\nabla u \cdot \nabla \upsilon\,dx\quad \mathrm{for\,\,all}\quad \upsilon \in V. 
\end{align*}
	$A$ is defined as 
	\begin{align*}
		Au:=-\mathcal{P}\Delta u,~u \in D(A)=V \cap H^2,
	\end{align*}
	where the operator $\mathcal{P}:L^2 \rightarrow H$ is the Leray orthogonal projection. 
	
	The nonlinear operator $B: V \rightarrow V'$ is defined as
	\begin{align*}
		B(u)=B(u,u)=\mathcal{P}[(u \cdot \nabla)u],
	\end{align*}
	and the trilinear form is defined as
	\begin{align*}
	b(u,\upsilon,w)=\left((u\cdot\nabla\upsilon),w\right)\quad \mathrm{for\,\,all}\quad u,\upsilon,w\in V.
    \end{align*}
	It is clear that the trilinear form $b(\cdot,\cdot,\cdot)$ satisfies the equalities
	\begin{align*}
		b(u,\upsilon,\upsilon)=0,~~~b(u,\upsilon,w)=-b(u,w,\upsilon),
	\end{align*}
	and the following estimates.
	\begin{lem}[See Proposition $9.2$ in \cite{R.T1986}]
		For $n=2$, if $u,~\upsilon,~w \in V$, then
		\begin{align*}
			|b(u,\upsilon,w)|\leq k\|u\|_{L^2}^{\frac{1}{2}}\|\nabla u\|_{L^2}^{\frac{1}{2}}\|\nabla \upsilon\|_{L^2}\|w\|_{L^2}^{\frac{1}{2}}\|\nabla w\|_{L^2}^{\frac{1}{2}}, 
		\end{align*}
		and if $u \in V$, $\upsilon \in D(A)$, and $w \in H$,then
		\begin{align*}
			|b(u,\upsilon,w)|\leq k\|u\|_{L^2}^{\frac{1}{2}}\|\nabla u\|_{L^2}^{\frac{1}{2}}\|\nabla \upsilon\|_{L^2}^{\frac{1}{2}}\|\Delta \upsilon\|_{L^2}^{\frac{1}{2}}\|w\|_{L^2}.
		\end{align*}
	\end{lem}
	Next, we recall  Gronwall's inequality, which provides exponentially decaying bounds. (See Lemma $2.8$ in \cite{R.T1986}). 
	\begin{lem}[\textbf{Gronwall's inequality}]\label{lem1}
		Let $t_0$ be a real number and $\varphi(t)\in L_{loc}^{1}(t_{0},\infty )$ be a function satisfying $\varphi'(t)\in L_{loc}^{1}(t_{0},\infty )$. Let $h(t)\in L_{loc}^{1}(t_{0},\infty)$ and $k(t)\in L_{loc}^{1}([t_{0},\infty ))$.  Let $\varphi(t)$ satisfy the differential inequality
		$$\displaystyle{ \varphi'(t) \leq \varphi(t)k(t) + h(t)\quad \mathrm{for\ a.e.}\quad t \in (t_{0},\infty ), }$$
		then
		$$\displaystyle{\varphi(t) \leq \varphi(t_{0})e^{\int_{t_{0}}^{t}k(\xi)\,d\xi } +\int_{t_{0}}^{t}h(\xi)e^{\int_{\xi}^{t}k(s)\,ds }\,d\xi\quad \mathrm{for\ all}\quad t \in [t_{0},\infty ). }$$
		In particular, if $k(t) = C$ for all $t\geq t_0$ is a real constant, then
		$$\displaystyle{ \varphi(t)\leq \varphi(t_{0})e^{-C(t_{0}-t) } + \int_{t_{0}}^{t}h(\xi)e^{-C(\xi-t)}\,d\xi\quad \mathrm{for\ all}\quad t \in [t_{ 0},\infty ). }$$
		Furthermore, if $h(t) = D$ for all $t\geq t_0$ is a real constant, then
		$$\displaystyle{ \varphi(t)\leq \varphi(t_{0})e^{-C(t_{0}-t) } + \frac{D}{C}(e^{-C(t_{0}-t) } - 1)\quad \mathrm{for\ all}\quad t \in [t_{ 0},\infty ). }$$
	\end{lem}
	Finally, we note that in this paper $C$ denotes a positive constant that may depend on the parameters of the problem or universal constants, $C_f$ denotes a positive constant depending on the parameter $f$, $R(\omega)$ denotes a random variable, and $R_{f,\nu}(\omega)$ denotes a random variable depending on $f$ and $\nu$. As for other constants, we denote them by the same pattern. Notice in particular that these constants may vary from line to line.

	\section{\large The Basic Set-up}\label{sec2}
	In this section, we first give some basic knowledge about random dynamical system (RDS). For a more detailed introduction about RDS, we can refer to \cite{I2002, H.F1994}.
	
	Let $X$ be a Banach space with Borel sigma-algebra and $(\Omega,\mathcal{F},P)$ be a probability space. If the following holds$:\theta_t:\Omega \to \Omega,~t \in \mathbb{R}$ is measurable, $\theta_0=id$(identity on $\Omega$), and $\theta_{t+s}=\theta_t \theta_s$ for all $t,s \in \mathbb{R}$, then $(\Omega,\mathcal{F},P,{\{\theta_t\}}_{t \in \mathbb{R}})$ is called a metric dynamical system. 
	
	Based on these, definitions of a continuous RDS, attraction, absorption and random attractors are given by the following four definitions, respectively.
	\begin{defi}[See Definition $2.1$ in \cite{H.F1994}]\label{def0}
		A continuous random dynamical system on X over a metric dynamical system $(\Omega,\mathcal{F},P,{\{\theta_t\}}_{t \in \mathbb{R}})$ is a measurable mapping 
		\begin{align*}
			\varphi:\mathbb{R}^{+} \times \Omega \times X \rightarrow X,~~(t,\omega,x) \mapsto \varphi(t,\omega,x),
		\end{align*}
		satisfying for $P$-a.e.$~\omega \in \Omega:$\par
		$(i)~\varphi(0,\omega)=id~\mathrm{on}~X;$\par
		$(ii)~\varphi(t+s,\omega)=\varphi(t,\theta_s\omega)  \varphi(s,\omega)$ for all $t,s\in \mathbb{R}^+;$\par
		$(iii)~\varphi(t,\omega):X \to X$ is continuous.\newline
		A family of maps satisfyng $(ii)$ ia called a cocycle.
	\end{defi}
	\begin{defi}[See Definition $3.3$ in \cite{H.F1994}]
		Given a random set $\mathscr{A}(\omega)$ and $B$ is another random set, we say $\mathscr{A}(\omega)$ attracts $B$ if $P$-a.s.
		\begin{align*}
			{\lim_{t\to + \infty}}dist (\varphi(t,\theta_{-t}\omega)B,\mathscr{A}(\omega))=0,
		\end{align*}
		where dist$(\cdot,\cdot)$ denotes the Hausdorff semi-distance. 
	\end{defi}
	\begin{defi}[See Definition $3.5$ in \cite{H.F1994}]
	Let $\mathscr{A}(\omega)$ and $B$ be random sets, for $P$-a.e.$~\omega \in \Omega$ if there exists $t_B(\omega)$ such that for all $t \geq  t_B(\omega)$,
	\begin{align*}
		\varphi(t,\theta_{-t}\omega)B \subset \mathscr{A}(\omega),
	\end{align*}
	then we say $\mathscr{A}(\omega)$ absorbs $B$.
	\end{defi}
	\begin{defi}[See Definition $3.9$ in \cite{H.F1994}]
		A random set $\mathscr{A}(\omega)$ is said to be a random attractor for the $\mathrm{RDS}$ $\varphi$ if \par
		$(i)\mathscr{A}(\omega)$ is a random compact set; \par
		$(ii)\mathscr{A}(\omega)$ is invariant;     \par
		$(iii)\mathscr{A}(\omega)$ attracts all deterministic bounded sets $B \subset X$.
	\end{defi}
	The following theorem is the main result in this paper, that is, the existence theorem of a random attractor for a continuous RDS.
	\begin{thm}[See Theorem $3.11$ in \cite{H.F1994}]\label{thm1}
		If there exists a random compact set absorbing every bounded set $B \subset X$, then the $\mathrm{RDS}$ possesses a random attractor $\mathcal A(\omega)$,
		$$\displaystyle{\mathcal A(\omega)=\overline{\mathop {\bigcup }\limits_{B \subset X} \Lambda_B(\omega)}},$$
		where $\Lambda_B(\omega)=\mathop {\bigcap }\limits_{s \geq 0} \overline{\mathop {\bigcup }\limits_{t \geq s}\varphi(t,\theta_{-t}\omega)B}$ is the omega-limit set of $B$.
	\end{thm}
	Next, applying Leray projection $\mathcal{P}$ onto the first equation in \eqref{a2} to get
	\begin{equation*}
		\frac{du}{dt} + B(u) + \nu Au  =-\rho \nabla \Phi +f+ \varepsilon u\circ\ \frac{dW(t)}{dt}. 
	\end{equation*}
	
	Finally, we generate the random dynamical system concerning SNPNS. Define $\{{\theta_t}\}_{t \in \mathbb{R}}$ by
	\begin{align*}
		\theta_t \omega(\cdot)=\omega(\cdot+t)-\omega(t),t \in \mathbb{R}, ~ \omega \in \varOmega.
	\end{align*}
	Hence, $(\Omega,\mathcal{F},P,{\{\theta_t\}}_{t \in \mathbb{R}})$ is a metric dynamical system. Moreover, there exists a $\theta_t$-invarient set $\tilde{\Omega \subseteq \Omega}$ of full $P$ measure such that for every $\omega$ $\in$ $\tilde{\Omega}$, the following holds:
	\begin{align}\label{l1}
		\frac{w(t)}{t}  \to  0~~\mathrm{as}~~  t  \to \pm \infty.
	\end{align}
	For convenience, we will not distinguish $\Omega$ and $\tilde{\Omega}$. Given $t \in \mathbb{R}$ and $\omega$ $\in$  $\Omega$, let the process
	\begin{align*}
		z(t,\omega)=e^{-\varepsilon\omega(t)}.
	\end{align*}
	Obviously, $z$ satisfies the equation
	\begin{align*}
		dz=-\varepsilon z \circ dW.
	\end{align*}
	Let $v^\varepsilon$ be a new variable given by 
	\begin{align*}
		v^\varepsilon=zu.
	\end{align*}
	Then, we obtain the new equations without stochastic differential
	\begin{equation}\label{b1}
		\left\{
		\begin{aligned} 
			\displaystyle{\frac{dv^\varepsilon}{dt} +  \nu\ Av^\varepsilon + \frac{1}{z(t,\omega)}B(v^\varepsilon)=-z(t,\omega)\rho\nabla \Phi+z(t,\omega)f},  \\
			\nabla \cdot v^\varepsilon = 0,  \\
			\rho=c_1-c_2, \\
			-\varepsilon_0 \Delta\Phi=\rho, \\
			\partial_t{c_1}+z^{-1}(v^\varepsilon\cdot\nabla c_1)=D\Delta c_1 +D\nabla\cdot(c_1\nabla \Phi), \\
			\partial_t{c_2}+z^{-1}(v^\varepsilon\cdot\nabla c_2)=D\Delta c_2-D\nabla\cdot(c_2\nabla \Phi). \\
		\end{aligned}
		\right.
	\end{equation}
	Let $\sigma=c_1+c_2$, system \eqref{b1} can be rewritten in its equivalent form as follows:
	\begin{equation}
		\left\{
		\begin{aligned}\label{b2} 
			\displaystyle{\frac{dv^\varepsilon}{dt} +  \nu\ Av^\varepsilon + \frac{1}{z(t,\omega)}B(v^\varepsilon)   =-z(t,\omega)\rho\nabla \Phi+z(t,\omega)f},  \\
			\nabla \cdot v^\varepsilon = 0,  \\
			\rho=c_1-c_2, \\
			-\varepsilon_0 \Delta\Phi=\rho, \\
			\partial_{t}\sigma+z^{-1}(v^\varepsilon\cdot\nabla \sigma)=D\Delta \sigma +D\nabla\cdot(\rho\nabla \Phi), \\
			\partial_{t}\rho+z^{-1}(v^\varepsilon\cdot\nabla \rho)=D\Delta \rho +D\nabla\cdot(\sigma\nabla \Phi). \\
		\end{aligned}
		\right.
	\end{equation}

	\section{\large Well-posedness of Solutions}\label{sec3}
	In this section, we will focus on the global well-posedness of solutions for system \eqref{b1}. The following calculations can be done rigorously by using a Faedo-Galerkin approximation method in \cite{R.T2012}.
	
	First of all, we give the existence and uniqueness theorem of local solutions with respect to the stochastic NPNS system.
	\begin{thm}[\textbf{Local Solution}]\label{thm2}	
		Suppose $v^\varepsilon(0) \in H^1$ and $c_i(0) \in {L}^{2}$. Then, there exists $T_0$ depending on $\|v^\varepsilon(0)\|_{H^1}$, $\|c_i(0)\|_{L^2}$, and the parameters of the problem such that system \eqref{b1} has a unique solution obeying $v^\varepsilon \in	L^{\infty}(0,T;H^1) \cap L^2(0,T;H^2)$ and $c_i \in L^{\infty}(0,T;L^2) \cap L^2(0,T;H^1)$ on $[0,T_0]$. 
		\begin{proof}
			Taking the ${L}^{2}$ inner product of the first equation in \eqref{b1} with $-\Delta v^\varepsilon $, we obtain
			\begin{align}\label{b3}
				\frac{1}{2}&\frac{d}{dt}\|\nabla v^\varepsilon(t)\|_{L^2}^{2}\ +\nu \|\Delta v^\varepsilon\|_{L^2}^{2}\notag \\&=\frac{1}{z(t,\omega)} \int (v^\varepsilon \cdot \nabla v^\varepsilon)\cdot \Delta v^\varepsilon+z(t,\omega)\int \rho\nabla \Phi\cdot\Delta v^\varepsilon-z(t,\omega) \int  f \Delta v^\varepsilon.
			\end{align}
			Using the fact that $v^\varepsilon$ is divergence free and integrating by parts, we get 
			\begin{align*}
				|\frac{1}{z(t,\omega)}\int (v^\varepsilon \cdot \nabla v^\varepsilon)\cdot \Delta v^\varepsilon| \leq \frac{1}{z(t,\omega)}\|\nabla v^\varepsilon \|_{L^4}^{2}\|\nabla v^\varepsilon \|_{L^2}.
			\end{align*}
			In view of Ladyzhenskaya's and Young's inequalities, we finally estimate
			\begin{align}\label{j1}
				\frac{1}{z(t,\omega)}\|\nabla v^\varepsilon \|_{L^4}^{2}\|\nabla v^\varepsilon \|_{L^2} &\leq  \frac{C}{z(t,\omega)}\|\nabla v^\varepsilon \|_{L^2}^{2}\|\Delta v^\varepsilon \|_{L^2}\notag \\ &\leq \frac{\nu}{6}\|\Delta v^\varepsilon \|_{L^2}^2+\frac{C}{z^2(t,\omega)}\|\nabla v^\varepsilon \|_{L^2}^{4}.
			\end{align}
			Thanks to elliptic regularity of the fourth equation in \eqref{b1}
			\begin{align*}
				\|\nabla \Phi\|_{L^\infty} \leq C\|\rho\|_{L^4},
			\end{align*}
			we estimate
			\begin{align*}
				\int \rho\nabla \Phi\cdot\Delta v^\varepsilon \leq \|\nabla \Phi\|_{L^\infty}\|\rho \|_{L^2}\|\Delta v^\varepsilon \|_{L^2} \leq C \|\rho\|_{L^4}\|\rho \|_{L^2}\|\Delta v^\varepsilon \|_{L^2}.
			\end{align*}
			Applying Ladyzhenskaya's inequality, we finally have
			\begin{align*}
				\int \rho\nabla \Phi\cdot\Delta v^\varepsilon \leq C \|\nabla \rho\|_{L^2}^{\frac{1}{2}}\|\rho \|_{L^2}^{\frac{3}{2}}\|\Delta v^\varepsilon \|_{L^2}.
			\end{align*}
			Using Young's inequality again yields
			\begin{align}\label{j2}
				z(t,\omega)\int \rho\nabla \Phi\cdot\Delta v^\varepsilon \leq \frac{\nu}{6}\|\Delta v^\varepsilon \|_{L^2}^2+C{z^2 (t,\omega)}\|\nabla \rho \|_{L^2}\|\rho \|_{L^2}^3.
			\end{align}
			Similarly, we also estimate
			\begin{align}\label{j3}
				|z(t,\omega) \int  f \Delta v^\varepsilon| \leq z(t,\omega)\|f\|_{L^2}\|\Delta v^\varepsilon\|_{L^2} \leq \frac{\nu}{6}\|\Delta v^\varepsilon \|_{L^2}^2+C_f+z^4(t,\omega).
			\end{align}
			Adding \eqref{j1}, \eqref{j2} and \eqref{j3}, we obtain
			\begin{align}\label{b4}
				\frac{d}{dt}&\|\nabla v^\varepsilon\|_{L^2}^{2}\ +\nu \|\Delta v^\varepsilon\|_{L^2}^{2}\notag \\ &\leq \frac{D}{2}\|\nabla \rho \|_{L^2}^2+C{z}^4(t,\omega)\|\rho \|_{L^2}^6+\frac{C}{z^2(t,\omega)}\|\nabla v^\varepsilon \|_{L^2}^{4}+C_f+Cz^4(t,\omega).
			\end{align}
			Let $\sigma=c_1+c_2$. Then, $\sigma$ and $\rho$ satisfy
			\begin{equation}
				\left\{
				\begin{aligned}\label{b5} 
					\partial_{t}\sigma+u\cdot\nabla \sigma=D\Delta \sigma +D\nabla\cdot(\rho\nabla \Phi), \\
					\partial_{t}\rho+u\cdot\nabla \rho=D\Delta \rho +D\nabla\cdot(\sigma\nabla \Phi). \\
				\end{aligned}
				\right.
			\end{equation}
			Taking the ${L}^{2}$ inner product of the first equation in \eqref{b5} with $\sigma$ and of the second equation with $\rho$, adding them, we obtain
			\begin{align}\label{ad1}
				\frac{1}{2}\frac{d}{dt}(\|\sigma\|_{L^2}^{2}\ + \|\rho\|_{L^2}^{2}\ ) =-D (\|\nabla\sigma\|_{L^2}^{2}\ + \|\nabla\rho\|_{L^2}^{2})+ \int \rho\Delta \Phi\sigma.
			\end{align}
			Using the fourth equation in \eqref{b1}
			\begin{align*}
				-\varepsilon_0 \Delta\Phi=\rho,
			\end{align*}
			we find
			\begin{align*}
				\int \rho\Delta \Phi\sigma \leq C\| \rho \|_{L^4} \|\rho\|_{L^2} \|\sigma\|_{L^4}.
			\end{align*}
			Thanks to Ladyzhenskaya's and Young's inequalities, we estimate 
			\begin{align}\label{ad2}
				\int \rho\Delta \Phi\sigma \leq \frac{D}{2}(\|\nabla\rho\|_{L^2}^{2} + \|\nabla\sigma\|_{L^2}^{2})+C\|\sigma\|_{L^2}^{4}+C\|\rho\|_{L^2}^{4}.
			\end{align}
			Then, adding \eqref{ad1} and \eqref{ad2}, we obtain
			\begin{align}\label{b6}
				\frac{d}{dt}( \|\sigma\|_{L^2}^{2}\ + \|\rho\|_{L^2}^{2}\ ) +D(\|\nabla\rho\|_{L^2}^{2} + \|\nabla\sigma\|_{L^2}^{2}) \leq C\|\sigma\|_{L^2}^{4}+C\|\rho\|_{L^2}^{4}.
			\end{align}
			Let 
			\begin{align*}
				M(t)=\|\nabla v^\varepsilon\|_{L^2}^{2}+\|\rho\|_{L^2}^{2}+\|\sigma\|_{L^2}^{2}.
			\end{align*}
			Combining \eqref{b4} and \eqref{b6} yields
			\begin{align}\label{b7}
				M'&(t)+\frac{D}{2}(\|\nabla\rho\|_{L^2}^{2} + \|\nabla\sigma\|_{L^2}^{2})+\nu \|\Delta v^\varepsilon\|_{L^2}^{2}\notag \\& \leq C{z}^8(t,\omega)+ \frac{C}{{z}^4(t,\omega)}+C{M(t)}^6+C_f+z^4(t,\omega).
			\end{align}
			This differential inequality guarantees the boundedness of the desired norms.
			Consequently, we conclude that 
			\begin{align*}
				v^\varepsilon \in L^{\infty}(0,T;H^1) \cap L^2(0,T;H^2)~~\mathrm{and}~~c_i \in L^{\infty}(0,T;L^2) \cap L^2(0,T;H^1) .
			\end{align*}
			Then we prove the uniqueness. For uniqueness, suppose $(v^\varepsilon_1,c_1^1,c_2^1)$ and $(v^\varepsilon_2,c_1^2,c_2^2)$ are two solutions of system \eqref{b1}. Let $\rho_1=c_1^1-c_2^1$, $\rho_2=c_1^2-c_2^2$, $\sigma_1=c_1^1+c_2^1$, $\sigma_2=c_1^2+c_2^2$. We write $v^\varepsilon=v^\varepsilon_1-v^\varepsilon_2$, $\rho=\rho_1-\rho_2$, $\sigma=\sigma_1-\sigma_2$. Then $v^\varepsilon$, $\sigma$, $\rho$ obey the system 
			\begin{equation}
				\left\{
				\begin{aligned}\label{b8}
					\displaystyle{\partial_{t}v^\varepsilon +  \nu\ Av^\varepsilon 
						=- \frac{1}{z(t,\omega)}\{B(v^\varepsilon_1)-B(v^\varepsilon_2)\}+ ({\rho}_2 \nabla {\Phi}_2-{\rho}_1\nabla {\Phi}_1 )z(t,\omega)}, \\
					\partial_{t}\sigma+\frac{1}{z(t,\omega)}(v^\varepsilon_1\cdot\nabla {\sigma}_1-v^\varepsilon_2\cdot\nabla {\sigma}_2)=D\Delta \sigma +D\nabla\cdot({\rho}_1 \nabla {\Phi}_1-{\rho}_2 \nabla {\Phi}_2), \\
					\partial_{t}\rho+\frac{1}{z(t,\omega)}(v^\varepsilon_1\cdot\nabla {\rho}_1-v^\varepsilon_2\cdot\nabla {\rho}_2)=D\Delta \rho +D\nabla\cdot({\sigma}_1\nabla {\Phi}_1-{\sigma}_2\nabla {\Phi}_2). \\
				\end{aligned}
				\right.
			\end{equation}
			We take the $L^2$ inner product of the first equation of \eqref{b8} with $v^\varepsilon$ and get 
			\begin{align}\label{b9}
				\frac{1}{2}&\frac{d}{dt} \|v^\varepsilon\|_{L^2}^{2}\ + \nu \|\nabla v^\varepsilon\|_{L^2}^{2}\notag\\& =-\frac{1}{z(t,\omega)}\langle B(v^\varepsilon_1)-B(v^\varepsilon_2),v^\varepsilon \rangle +z(t,\omega)\langle {\rho}_2\nabla {\Phi}_2-{\rho}_1\nabla {\Phi}_1,v^\varepsilon \rangle.
			\end{align}
			We estimate the term 
			\begin{align*}
				|\frac{1}{z(t,\omega)}\langle B(v^\varepsilon_1)-B(v^\varepsilon_2),v^\varepsilon \rangle|&=|\frac{1}{z(t,\omega)} \int (v^\varepsilon_1\cdot\nabla v^\varepsilon_1-v^\varepsilon_2\cdot\nabla v^\varepsilon_2)\cdot v^\varepsilon dx|\\
				&=|\frac{1}{z(t,\omega)} \int (v^\varepsilon\cdot\nabla v^\varepsilon_1+v^\varepsilon_2\cdot\nabla v^\varepsilon)\cdot v^\varepsilon  dx|.  
			\end{align*}
			Applying Ladyzhenskaya's and Young's inequalities, we have
			\begin{align}\label{b10}
				|\frac{1}{z(t,\omega)} \int (v^\varepsilon\cdot\nabla v^\varepsilon_1+v^\varepsilon_2\cdot\nabla v^\varepsilon)\cdot v^\varepsilon dx|  &\leq \frac{C}{z(t,\omega)} \|v^\varepsilon\|_{L^2}^{\frac{3}{2}} \|\nabla v^\varepsilon\|_{L^2}^{\frac{1}{2}}\|\nabla v^\varepsilon_1\|_{L^2}^{\frac{1}{2}}\|\Delta v^\varepsilon_1\|_{L^2}^{\frac{1}{2}}\notag\\&\leq \frac{\nu}{4} \|\nabla v^\varepsilon\|_{L^2}^{2}+{\frac{C}{z^{\frac{4}{3}}(t,\omega)}}\|v^\varepsilon\|_{L^2}^2\|\nabla v^\varepsilon_1\|_{L^2}^{\frac{2}{3}}\|\Delta v^\varepsilon_1\|_{L^2}^{\frac{2}{3}}. 
			\end{align}
			We estimate the term 
			\begin{align}\label{b11}
				|z(t,\omega)\langle {\rho}_2\nabla {\Phi}_2-{\rho}_1\nabla {\Phi}_1,v^\varepsilon \rangle|=|z(t,\omega)\int ({\rho}\nabla {\Phi_1}+{\rho}_2\nabla {\Phi})\cdot v^\varepsilon dx|.
			\end{align}
			In view of elliptic regularity
			\begin{align}\label{b12}
				\|\nabla \Phi\|_{L^\infty} \leq C\|\rho\|_{L^4},
			\end{align}
			we get
			\begin{align}\label{b13}
				|z&(t,\omega)\int ({\rho}\nabla {\Phi_1}+{\rho}_2\nabla {\Phi})\cdot v^\varepsilon dx|\notag\\&  \leq Cz(t,\omega)[\|\nabla \Phi_1\|_{L^\infty}\|\rho\|_{L^2}\|v^\varepsilon\|_{L^2}+\|\rho_2\|_{L^2}\|\rho\|_{L^2}^{\frac{1}{2}}\|\nabla \rho\|_{L^2}^{\frac{1}{2}}\|v^\varepsilon\|_{L^2}].		
			\end{align}
			Then, we take the $L^2$ inner product of the second equation of \eqref{b8} with $\sigma$ and obtain
			\begin{align}\label{b14}
				\frac{1}{2}&\frac{d}{dt} \|\sigma\|_{L^2}^{2}\ + D \|\nabla \sigma\|_{L^2}^{2}\notag \\& =-\frac{1}{z(t,\omega)}\int (v^\varepsilon_1\cdot\nabla \sigma_1-v_2\cdot\nabla \sigma_2) \sigma+D\int (\nabla \cdot (\rho_1\nabla \Phi_1-\rho_2\nabla \Phi_2))\sigma.
			\end{align}
			We have 
			\begin{align}\label{b15}
				|\frac{1}{z(t,\omega)}\int (v^\varepsilon_1\cdot\nabla \sigma_1-v^\varepsilon_2\cdot\nabla \sigma_2) \sigma|&=\frac{1}{z(t,\omega)}|\int (v^\varepsilon\cdot\nabla \sigma_1+v^\varepsilon_2\cdot\nabla \sigma) \sigma|\notag \\&\leq \frac{C}{z(t,\omega)}\|\nabla \sigma_1\|_{L^2}\|v^\varepsilon\|_{L^2}^{\frac{1}{2}}\|\nabla v^\varepsilon\|_{L^2}^{\frac{1}{2}}\|\sigma\|_{L^2}^{\frac{1}{2}}\|\nabla \sigma\|_{L^2}^{\frac{1}{2}},
			\end{align}
			and
			\begin{align}\label{b16}
				|\int (\nabla \cdot (\rho_1\nabla \Phi_1-\rho_2\nabla \Phi_2))\sigma|&=|\int (\nabla \cdot (\rho\nabla \Phi_1+\rho_2\nabla \Phi))\sigma|\notag\\ &\leq C[\|\nabla \Phi_1\|_{L^\infty}\|\rho \|_{L^2}\|\nabla \sigma\|_{L^2}+\|\rho_2 \|_{L^2}\|\nabla \Phi\|_{L^\infty}\|\nabla \sigma \|_{L^2}]\notag\\&\leq C[\|\nabla \Phi_1\|_{L^\infty}\|\rho \|_{L^2}\|\nabla \sigma\|_{L^2}+\|\rho_2 \|_{L^2}\|\rho \|_{L^2}^{\frac{1}{2}}\|\nabla \rho \|_{L^2}^{\frac{1}{2}}\|\nabla \sigma \|_{L^2}].
			\end{align}
			Finally, we take the $L^2$ inner product of the third equation of \eqref{b8} with $\rho$ and obtain
			\begin{align}\label{b17}
				\frac{1}{2}&\frac{d}{dt} \|\rho\|_{L^2}^{2}\ + D \|\nabla \rho\|_{L^2}^{2}\notag\\& =-\frac{1}{z(t,\omega)}\int (v^\varepsilon_1\cdot\nabla \rho_1-v^\varepsilon_2\cdot\nabla \rho_2)\rho+D\int (\nabla \cdot (\sigma_1\nabla \Phi_1-\sigma_2\nabla \Phi_2))\rho.
			\end{align}
			Estimating the terms as above, we also have
			\begin{align}\label{b18}
				|\frac{1}{z(t,\omega)}\int (v^\varepsilon_1\cdot\nabla \rho_1-v^\varepsilon_2\cdot\nabla \rho_2)\cdot \rho|&=\frac{1}{z(t,\omega)}|\int (v^\varepsilon\cdot\nabla \rho_1+v^\varepsilon_2\cdot\nabla \rho)\cdot \rho|\notag\\ &\leq \frac{C}{z(t,\omega)}\|\nabla \rho_1\|_{L^2}\|v^\varepsilon\|_{L^2}^{\frac{1}{2}}\|\nabla v^\varepsilon\|_{L^2}^{\frac{1}{2}}\|\rho\|_{L^2}^{\frac{1}{2}}\|\nabla \rho\|_{L^2}^{\frac{1}{2}},
			\end{align}
			and
			\begin{align}\label{b19}
				|\int (\nabla \cdot (\sigma_1\nabla \Phi_1-\sigma_2\nabla \Phi_2))\rho|&=|\int (\nabla \cdot (\sigma\nabla \Phi_1+\sigma_2\nabla \Phi))\rho| \notag\\&\leq C[\|\nabla \Phi_1\|_{L^\infty}\|\sigma \|_{L^2}\|\nabla \rho\|_{L^2}+\|\sigma_2 \|_{L^2}\|\rho \|_{L^2}^{\frac{1}{2}}\|\nabla \rho \|_{L^2}^{\frac{3}{2}}].
			\end{align}
			Adding \eqref{b10} to \eqref{b19} together and applying Young's inequality, we deduce
			\begin{align}\label{b20}
				\frac{d}{dt}[\|v^\varepsilon\|_{L^2}^2+\|\rho\|_{L^2}^2+\|\sigma\|_{L^2}^2] \leq CC(t)[\|v^\varepsilon\|_{L^2}^2+\|\rho\|_{L^2}^2+\|\sigma\|_{L^2}^2],
			\end{align}
			where
			\begin{align*}
				C(t)=&\|\nabla v^\varepsilon_1\|_{L^2}\|\Delta v^\varepsilon_1\|_{L^2}+\|\nabla \rho_1\|_{L^2}^{4}+\|\nabla \sigma_1\|_{L^2}^{4}+\| \sigma_2\|_{L^2}^{4}+\|\rho _2\|_{L^2}^{4}\notag \\&+{z^2(t,\omega)}+\frac{1}{{z^4(t,\omega)}}+1.
			\end{align*}
			Since 
			\begin{align*}
				\int_{0}^{t} C(s)ds \textless \infty,
			\end{align*}
			for any $t \in[0,T_0]$, then we can prove the uniqueness.
		\end{proof}
	\end{thm} 
	Next, the following theorem shows that the solution is continuous corresponding to initial data.
	\begin{thm}\label{thm3}
		Suppose $v^\varepsilon(0) \in H^1$ and $c_i(0) \in {L}^{2}$. Then, the solution of \eqref{b1} is continuous in initial data.
		\begin{proof}
			Applying Lemma \ref{lem1} to \eqref{b20}, we get
			\begin{align}\label{b21}
				\|v^\varepsilon\|_{L^2}^2+\|\rho\|_{L^2}^2+\|\sigma\|_{L^2}^2 \leq e^{C\int_{0}^{T_0} C(s)ds}(\|v^\varepsilon(0)\|_{L^2}^2+\|\rho(0)\|_{L^2}^2+\|\sigma(0)\|_{L^2}^2).
			\end{align}
			Hence the proof is completed.
		\end{proof}
	\end{thm}
	Then, we explore the existence of global regular solutions under the following condition \eqref{b22}.
	\begin{thm}\label{thm4}
		Suppose $v^\varepsilon(0) \in H^1$ and $c_i(0) \in H^1$. Let $T > 0$. Suppose $(v^\varepsilon,c_1,c_2)$ is the solution of \eqref{b1} on the interval $[0,T]$ with 
		\begin{align}\label{b22}
			\int_{0}^{T}(\|c_1(t)\|_{L^2}^2+\|c_2(t)\|_{L^2}^2)dt \textless \infty.
		\end{align}
		Then, $v^\varepsilon \in 
		L^{\infty}(0,T;H^1) \cap L^2(0,T;H^2)$ and $c_i \in L^{\infty}(0,T;H^1) \cap L^2(0,T;H^2) $. 
		\begin{proof}
			The differential inequality \eqref{b6} gives
			\begin{align}\label{b23}
				\frac{d}{dt}( \|\sigma\|_{L^2}^{2}\ + \|\rho\|_{L^2}^{2}\ ) +D(\|\nabla\rho\|_{L^2}^{2} + \|\nabla\sigma\|_{L^2}^{2}) \leq C\|\sigma\|_{L^2}^{4}+C\|\rho\|_{L^2}^{4} \leq C(\|\sigma\|_{L^2}^{2}+C\|\rho\|_{L^2}^{2})^2.
			\end{align}
			Under the assumption \eqref{b22}, we conclude that $c_i \in L^{\infty}(0,T;L^2) \cap L^2(0,T;H^1) $.  \newline
			The differential inequality \eqref{b4} gives
			\begin{align}\label{b24}
				\frac{d}{dt}&\|\nabla v^\varepsilon\|_{L^2}^{2}\ +\nu \|\Delta v^\varepsilon\|_{L^2}^{2}\notag \\&\leq \frac{D}{2}\|\nabla \rho \|_{L^2}^2+C{z}^4(t,\omega)\|\rho \|_{L^2}^6+\frac{C}{z^2(t,\omega)}\|\nabla v^\varepsilon \|_{L^2}^{4}+C_f+Cz^4(t,\omega).
			\end{align}
			Thus, we obtain that $v^\varepsilon \in  L^{\infty}(0,T;H^1) \cap L^2(0,T;H^2)$. \newline
			Taking the $L^2$ inner product of the equation satisfied by $\sigma$ in \eqref{b5} with -$\Delta \sigma$, we get
			\begin{align}\label{b25}
				\frac{1}{2}\frac{d}{dt} \|\nabla \sigma\|_{L^2}^{2}\ + D \|\Delta \sigma\|_{L^2}^{2} =\frac{1}{z(t,\omega)}\int (v^\varepsilon\cdot \nabla \sigma)\Delta \sigma -D\int \nabla\cdot(\rho \nabla \Phi)\Delta \sigma.
			\end{align}
			We estimate
			\begin{align}\label{b26}
				|\int \rho \Delta \Phi \Delta \sigma| \leq \frac{1}{6}\|\Delta \sigma\|_{L^2}^2+C\|\rho\|_{L^2}^4+C\|\nabla \rho\|_{L^2}^4,
			\end{align}
			\begin{align}\label{b27}
				|\int (\nabla \rho \cdot  \nabla \Phi) \Delta \sigma| \leq \frac{1}{6}\|\Delta \sigma\|_{L^2}^2+C\|\nabla \rho\|_{L^2}^4,
			\end{align}
			and 
			\begin{align}\label{b28}
				|\frac{1}{z(t,\omega)}\int (v^\varepsilon\cdot \nabla \sigma)\Delta \sigma| \leq \frac{D}{6}\|\Delta \sigma\|_{L^2}^2+\frac{C}{z^4(t,\omega)}+C\|\nabla v^\varepsilon\|_{L^2}^4\|\nabla \sigma\|_{L^2}^4,
			\end{align}
			where we used elliptic regularity, Young's inequality, Ladyzhenskaya's inequality and Poincar$\acute{\mathrm{e}}$'s inequality.
			
			Then, we take the $L^2$ inner product of the equation satisfied by $\rho$ in \eqref{b5} with -$\Delta \rho$ to obtain
			\begin{align}\label{b29}
				\frac{1}{2}\frac{d}{dt} \|\nabla \rho\|_{L^2}^{2}\ + D \|\Delta \rho\|_{L^2}^{2} =\frac{1}{z(t,\omega)}\int (v^\varepsilon\cdot \nabla \rho)\Delta \rho -D\int \nabla\cdot(\sigma \nabla \Phi)\Delta \rho.
			\end{align}
			Similarly, we estimate
			\begin{align}\label{b30}
				|\int \sigma \Delta \Phi \Delta \rho| \leq \frac{1}{6}\|\Delta \rho\|_{L^2}^2+C\|\sigma\|_{L^2}^2\|\nabla \sigma\|_{L^2}^2+C\|\nabla \rho\|_{L^2}^4,
			\end{align}
			\begin{align}\label{b31}
				|\int (\nabla \sigma \cdot  \nabla \Phi) \Delta \rho| \leq \frac{1}{6}\|\Delta \rho\|_{L^2}^2+C\|\nabla \rho\|_{L^2}^4+C\|\nabla \sigma\|_{L^2}^4,
			\end{align}
			and 
			\begin{align}\label{b32}
				|\frac{1}{z(t,\omega)}\int (v^\varepsilon\cdot \nabla \rho)\Delta \rho|&=|z^{-1}\int (\nabla v^\varepsilon \nabla \rho \nabla \rho|\notag\\&  \leq \frac{D}{6}\|\Delta \rho\|_{L^2}^2+\frac{C}{z^4(t,\omega)}+C\|\nabla v^\varepsilon\|_{L^2}^4\|\nabla \rho\|_{L^2}^4.
			\end{align}
			Adding \eqref{b25} to \eqref{b32} together, we conclude that $c_i \in L^{\infty}(0,T;H^1) \cap L^2(0,T;H^2) $ with bounds depending on the initial data and $T$.
		\end{proof}
	\end{thm}
    The following remark is cited to provide the needed estimates. 
	\begin{rem}[See \cite{P.M2019}]\label{Remark 1}
		Under the conditions of Theorem \ref{thm4}, if $c_i(0) \geq 0$, then $c_i(t) \geq 0$ for $0 \leq t \leq T.$
	\end{rem}
	Before presenting the final result in this section, we also need to give a priori $L^2$ bounds.
	\begin{prop}\label{Prop5}
		Let $v^\varepsilon(0) \in H^1$ and $c_i(0) \in H^1$. Assume that $c_i(t) \geq  0$ holds for all $t \geq 0 $. Then, there exists an absolute constant $C > 0 $ such that 
		\begin{align}\label{b33}
			&\int_{t}^{t+T}(\|\nabla \rho (s)\|_{L^2}^2+\|\nabla \sigma (s)\|_{L^2}^2+\frac{1}{\varepsilon}\|\nabla \rho (s)\|_{L^3}^3)ds\notag\\& \leq \frac{1}{2D}(2\|\sigma(t_0)\|_{L^2}^2+2\|\overline {\sigma}\|_{L^2}^2+\|\rho(t_0)\|_{L^2}^2)Te^{-2CDt},
		\end{align}
		for all $t \geq 0$, $T>0$. 
		\begin{proof}
			The proof is similar to Proposition $3$ in \cite{E.M2021}, so we omit it.
		\end{proof}
	\end{prop}
	Finally, with the results we have just proved, we are now ready to show that this local solution can be extended to a strong analytic solution on $[0,T]$ for any $T > 0$.
	\begin{thm}\label{thm5}
		Let $v^\varepsilon(0)$ $\in$  $H^1$ and $c_i(0) \in H^1$ be nonnegative with $c_i(0) \geq 0 $. Let $T > 0$. Then, there exists a unique solution $(v^\varepsilon,c_1,c_2)$ satisfying $v^\varepsilon \in 
		L^{\infty}(0,T;H^1) \cap L^2(0,T;H^2)$ and $c_i \in L^{\infty}(0,T;H^1) \cap L^2(0,T;H^2) $. Moreover $c_i(t) \geq 0$ holds on $[0,T]$. 
		\begin{proof}
			 By the local existence theorem (Theorem \ref{thm2}), there exists $T_0 \geq 0$ depending only on the norms of initial data in $H^1$ such that the solution exists and belongs to $H^1$. By Remark \ref{Remark 1}, we find $c_i(t) \geq 0$. Moreover, the inequality \eqref{b33} holds on $[0,T_0]$. It follows from Theorem \ref{thm4} that the solution is bounded in $H^1$. Applying the local existence theorem again and starting from $T_0$, we deduce that the solution can be extended for $T_1$ $\geq T_0$. The inequality \eqref{b33} holds on $[0,T_1]$. Since the inequality holds as long as $c_i$ $\geq 0$, reasoning by contradiction we conclude that the solution can extend to the whole interval $[0,T]$.
		\end{proof}
	\end{thm}

	\section{\large Random Attractors}\label{sec4}
	In the previous section, we have proved existence, uniqueness, and continuous dependence of solutions with respect to initial data. Let us denote the solution of system \eqref{b2} by
	\begin{align*}
	 (v^\varepsilon(t,\omega,t_0,v^\varepsilon_0),\sigma(t,\omega,t_0,\sigma_0),\rho(t,\omega,t_0,\rho_0)).
	\end{align*}
	By Theorem \ref{thm3}, $v^\varepsilon$ is continuous with respect to initial data and measurable in $\omega \in \Omega$. Thus it is clear that the mapping $(v^\varepsilon_0,\sigma_0,\rho_0)\mapsto(v^\varepsilon(t,\omega,t_0,v^\varepsilon_0),\sigma(t,\omega,t_0,\sigma_0),\rho(t,\omega,t_0,\rho_0))$ is continuous for all $t \geq t_0$. These guarantee that we can define a cocycle $S(t,\omega)$ by 
	\begin{align*}
		S(t,\omega)x_0= (u,\sigma,\rho)=(z^{-1}(t,\omega)v^\varepsilon(t,\omega,0,u_0),\sigma(t,\omega,0,\sigma_0),\rho(t,\omega,0,\rho_0)),~t \geq 0,
	\end{align*}
	where $x_0= ( u_0,\sigma_0,\rho_0 )$.
	By Theorem \ref{def0}, it is easy to find that $S(t,\omega)$ is a continuous RDS on $H$ over ${\{\theta_t\}}_{t \in \mathbb{R}}$ and $(\Omega,\mathcal{F},P,{\{\theta_t\}}_{t \in \mathbb{R}})$, which implies that 
	\begin{align*}
		S(t,\theta_{-t}\omega) x_0=S(0,\omega,-t,x_0).
	\end{align*}
		
	In this section, we proceed in several steps. First it can be shown that the random dynamical system associated with system \eqref{b2} has a random absorbing set in $H$. Next we are able to prove that there also exists a random absorbing set in $V$. Finally, the existence of a compact random attractor in $H$ is a direct consequence by making use of the Sobolev embedding theorem and Theorem \ref{thm1}. For more basic theory concerning random attractors, we can refer to \cite{Z.M.F1993, H.F1994, H.A.F1997, F.B1996,R.Z2017}. In order to prove these conclusions, we need to provide some auxiliary lemmas as follows.
	\begin{lem}\label{lem5}
		Let $(v^\varepsilon,\sigma,\rho)$ be a solution of system \eqref{b2}. Let $v^\varepsilon(t_0)$ $\in$  $H^1$ and $c_i(t_0) \in {L}^{2}$. Then,
		\begin{align}\label{d1}
			\|\sigma\|_{L^2}^{2}\ + \|\rho\|_{L^2}^{2}\ \leq (\|\sigma(t_0)\|_{L^2}^{2}\ + \|\rho(t_0)\|_{L^2}^{2})e^{2D(t_0-t)},~t \geq t_0.
		\end{align}
		\begin{proof}
	We recall $\sigma$  and $\rho$ obey
			\begin{equation}
				\left\{
				\begin{aligned}\label{d2} 
					\partial_{t}\sigma+z^{-1}(v^\varepsilon\cdot\nabla \sigma)=D\Delta \sigma +D\nabla\cdot(\rho\nabla \Phi) \\
					\partial_{t}\rho+z^{-1}(v^\varepsilon\cdot\nabla \rho)=D\Delta \rho +D\nabla\cdot(\sigma\nabla \Phi) \\
				\end{aligned}
				\right.
			\end{equation}
			Taking the ${L}^{2}$ inner product of the first equation in \eqref{d2} with $\sigma$ and of the second equation with $\rho$, we add them and have
			\begin{align*}
				\frac{1}{2}\frac{d}{dt}( \|\sigma\|_{L^2}^{2}\ + \|\rho\|_{L^2}^{2}\ ) =-D (\|\nabla\sigma\|_{L^2}^{2}\ + \|\nabla\rho\|_{L^2}^{2}\ )+ \int \rho\Delta \Phi\sigma.
			\end{align*}
			In view of the fourth equation in \eqref{b2}, we know that
			\begin{align*}
				\int \rho\Delta \Phi\sigma=-\frac{1}{\varepsilon_0} \int {\rho}^2\sigma,
			\end{align*}
			where $\sigma > 0$. Then, we get
			\begin{align}\label{d3}
				\frac{d}{dt}( \|\sigma\|_{L^2}^{2}\ + \|\rho\|_{L^2}^{2}\ ) +2D(\|\nabla\rho\|_{L^2}^{2} + \|\nabla\sigma\|_{L^2}^{2}) \leq 0.
			\end{align}
			Thanks to Lemma \ref{lem1}, we obtain \eqref{d1}.
		\end{proof}
	\end{lem}
	\begin{lem}\label{lem6}
		Let $v^\varepsilon(t_0)$ $\in$  $H^1$ and $c_i(t_0) \in H^1$. Then, there exists a positive constant $a$ depending on $D$, $\varepsilon_0$ and $\nu$, and a positive constant $A$ depending on the initial data, universal constants, $\varepsilon$ and $\omega$ such that 
		\begin{align}\label{d4}
			\|\nabla \rho (t)\|_{L^2}^2+\|\nabla \sigma (t)\|_{L^2}^2 \leq  Ae^{-a(t-t_0)},
		\end{align}
		for all $t \geq t_0$. 
		\begin{proof}
			Adding \eqref{b26} to \eqref{b32} together, the conclusion can be proved.
		\end{proof}
	\end{lem}
	Lemma \ref{lem5} and Lemma \ref{lem6} show that $\sigma$ and $\rho$ decay exponentially in $H$ and $V$.
	\begin{lem}\label{lem7}
		Let $(v^\varepsilon,\sigma,\rho)$ be a solution of system \eqref{b2}. Let $v^\varepsilon(t_0)$ $\in$  $H^1$ and $c_i(t_0) \in H^1$. Then,
		\begin{align}\label{d5}
			\|v^\varepsilon(t)\|_{L^2}^{2}\ \leq {e}^{-\nu t}[	\|v^\varepsilon(t_0)\|_{L^2}^{2}\ {e}^{\nu t_0} +\frac{2}{\nu}\int_{t_{0}}^{t} {z^{2}(s,w)}(\|\rho(s) \|_{L^2}^{3}\|\nabla \rho(s) \|_{L^2}+\|f\|_{L^2}^2) {e}^{\nu s}ds ],
		\end{align}
		for $t \geq t_0$.
		\begin{proof}
			Taking the ${L}^{2}$ inner product of the first equation in \eqref{b2} with $v^\varepsilon$ yields
			\begin{align}\label{jia1}
				\frac{1}{2}\|v^\varepsilon\|_{L^2}^{2}\ +\nu \|\nabla v^\varepsilon\|_{L^2}^{2}=-(\rho\nabla\Phi,v^\varepsilon)z(t,w)+(f,v^\varepsilon)z(t,w).
			\end{align}
			In view of Young's and Ladyzhenskaya's inequalities, we estimate
			\begin{align}\label{jia2}
				|(\rho\nabla\Phi,v^\varepsilon)z(t,w)| \leq \frac{\nu}{4}\|\nabla v^\varepsilon\|_{L^2}^{2}+\frac{1}{\nu}{z^{2}(t,w)}\|\rho \|_{L^2}^{3}\|\nabla \rho \|_{L^2},
			\end{align}
			and
			\begin{align}\label{jia3}
				|(f,v)z(t,w)| \leq  \frac{\nu}{4}\|\nabla v^\varepsilon\|_{L^2}^{2}+\frac{1}{\nu}{z^{2}(t,w)}\|f\|_{L^2}^2.
			\end{align}
			Putting \eqref{jia1} to \eqref{jia3} together, we get
			\begin{align}\label{d6}
				\frac{d}{dt} \|v^\varepsilon\|_{L^2}^{2}\ + \nu \|\nabla v^\varepsilon\|_{L^2}^{2} \leq\frac{2}{\nu}{z^{2}(t,w)}\|\rho \|_{L^2}^{3}\|\nabla \rho \|_{L^2}+\frac{2}{\nu}{z^{2}(t,w)}\|f\|_{L^2}^2.
			\end{align}
			By Lemma \ref{lem1}, we obtain \eqref{d5}.
		\end{proof}
	\end{lem}
	Lemma \ref{lem7} says $v^\varepsilon$ decays exponentially in $H$. By making full use of the lemmas above, we prove further that there exists a random absorbing set in $H$.
	\begin{lem}\label{lem8}
		There exist a positive number $R_0$ and a random variable $R_1(\omega)$ depending on the initial data, $\varepsilon_0$, $\nu$, $f$ and universal constants, such that: For every $E >0 $, there exists $t(\omega,E) \leq -1$, such that for all $x_0=(u_0,\rho_0,\sigma_0) \in {H^1}$ with $|x_0| < E$, and for any $t_0 < t(\omega,E)$, we have
		\begin{align}\label{d7}
			\|\sigma(t,\omega,t_0,\sigma_0)\|_{L^2} \leq R_0,~for~all~t \in [-1,0],
		\end{align}
		\begin{align}\label{d8}
			\|\rho(t,\omega,t_0,\rho_0)\|_{L^2} \leq R_0,~for~all~t \in [-1,0],
		\end{align}
		\begin{align}\label{d9}
			\|v^\varepsilon(t,\omega,t_0,z(t_0,\omega)u_0)\|_{L^2} \leq  R_1(\omega),~for~all~t \in [-1,0].
		\end{align}
		Therefore, $B(0,2R_0+R_1(\omega))$ is a random absorbing set in H.
		\begin{proof}
			By Lemma \ref{lem5}, for all $t \in [-1,0]$ and $t_0 \leq 1$, we get 
			\begin{align*}
				\|\sigma\|_{L^2}^{2}\ + \|\rho\|_{L^2}^{2}\ \leq (\|\sigma(t_0)\|_{L^2}^{2}\ + \|\rho(t_0)\|_{L^2}^{2}\ )e^{2D(t_0+1)} \leq 2{E^2}e^{2D(t_0+1)}.
			\end{align*}
			Choose $t_1$$(E)$ $\leq  -1$ such that 
			\begin{align*}
				2{E^2}e^{2D(t_1+1)} \leq {R_0}^2,
			\end{align*}
			provided $t_0 \leq t_1(E)$, which proves \eqref{d7} and \eqref{d8}.\newline
			We see from Lemma \ref{lem7} that 
			\begin{align}\label{d10}
				\|v^\varepsilon(t)\|_{L^2}^{2}\ \leq {e}^{-\nu t}[	{z}^{2}(t_0,\omega)	\|u_0\|_{L^2}^{2}\ {e}^{\nu t_0} +\frac{2}{\nu}\int_{t_{0}}^{t} {z^{2}(s,w)}(\|\rho(s) \|_{L^2}^{3}\|\nabla \rho(s) \|_{L^2}+\|f\|_{L^2}^2) {e}^{\nu s}ds],
			\end{align}
			for all $t \in [-1,0]$ and $t_0 \leq -1$. 
			Then by \eqref{l1}, we notice that
			\begin{align*}
				{z}^{2}(t_0,\omega)e^{\nu t_0} \to 0,~ t \to -\infty, 
			\end{align*}
			which implies that there exists $t_2(\omega,E) \leq -1$ such that if $t_0 \leq  t_2$,
			\begin{align*}
				{z}^{2}(t_0,\omega)	\|u_0\|_{L^2}^{2}\ {e}^{\nu t_0} \leq E^2{z}^{2}(t_0,\omega) {e}^{\nu t_0} \leq 1.
			\end{align*}
			Then, we estimate the integrals on the right hand of \eqref{d10}. \newline
			Let 
			\begin{align*}
				I&=\int_{t_{0}}^{t} {z(s,w)}^{2}\|\rho(s) \|_{L^2}^{3}\|\nabla \rho(s) \|_{L^2} {e}^{\nu s}ds\\ 
				&=(\int_{\frac{t_0}{2}}^{t} + \int_{t_0}^{\frac{t_0}{2}}) [{z^{2}(s,w)}\|\rho(s) \|_{L^2}^{3}\|\nabla \rho(s) \|_{L^2} {e}^{\nu s}]ds \\
				&=I_1+I_2.
			\end{align*}
			To estimate $I_1$, if $\frac{t_0}{2} \textless s \textless t$, then
			\begin{align*}
				\|\rho(s)\|_{L^2} \leq R_0,
			\end{align*}
			provided $t_0$ $\leq$ $t_3(E)$ for some chosen $t_3(E)$ $\leq -1$. \newline
			We see from Lemma \ref{lem6} that 
			\begin{align*}
				\|\nabla \rho (t)\|_{L^2}^2+\|\nabla \sigma (t)\|_{L^2}^2 \leq  Ae^{-a(t-t_0) }.
			\end{align*}
			If $\frac{t_0}{2}< s < t$, then
			\begin{align*}
				\|\nabla \rho (s)\|_{L^2}^2 \leq A,
			\end{align*}
			$A$ is given in lemma \ref{lem6}, provided $t_0$ $\leq$ $t_3(E)$ for some chosen $t_3(E)$ $\leq -1$. \newline			 
			For $t_0 \leq t_3(E)$, $t \in  [-1,0]$,
			\begin{align*}
				I_1 \leq \sqrt A{R_0}^3\int_{-\infty}^{0} {z}^{2}(s,w) {e}^{\nu s} ds = r_1(\omega),
			\end{align*}
			$r_1(\omega)$ is a random variable being independent of $E$. Since $\frac{w(t)}{t}$ $\to 0$ as $t \to$ $\pm$ $\infty$, the ergodic property implies that $s \mapsto$ ${z}^{2}(s,w) {e}^{\nu s}$ is pathwise integrable on $(-\infty,0]$.\newline
			To estimate $I_2$, we choose $t_4(E) \leq -1$ such that 
			\begin{align*}
				8\sqrt A E^3e^{\frac{\nu}{2}t_4} \leq 1.
			\end{align*}
			By Lemma \ref{lem5}, we deduce
			\begin{equation*}
				\|\rho(s)\|_{L^2}\ \leq 2E, 
			\end{equation*}
			for all s $\in [t_0,\frac{t_0}{2}]$. If $\frac{t_0}{2} \leq t_4(E)$, then we derive
			\begin{align*}
				I_2&=\int_{t_0}^{\frac{t_0}{2}} {z}^{2}(s,w)\|\rho(s)\|_{L^2}^{3}\|\nabla \rho(s)\|_{L^2} {e}^{\nu s}ds \\	&\leq \int_{t_0}^{\frac{t_0}{2}} [8\sqrt A E^3{e}^{\frac{\nu}{2} s}]{z}^{2}(s,w) {e}^{\frac{\nu}{2} s}ds \\
				&\leq \int_{t_0}^{\frac{t_0}{2}} {z}^{2}(s,w) {e}^{\frac{\nu}{2} s}ds \\
				&\leq \int_{-\infty}^{0} {z}^{2}(s,w) {e}^{\frac{\nu}{2} s}ds=r_2(\omega). 
			\end{align*}
			Let
			\begin{align*}
				J=\int_{t_{0}}^{t} {z(s,w)}^{2}\|f\|_{L^2}^{2}{e}^{\nu s}ds.
			\end{align*}
			Then, we estimate
			\begin{align*}
				J \leq \int_{-\infty}^{0} {z(s,w)}^{2}\|f\|_{L^2}^{2}{e}^{\nu s}ds=R_{f,\nu}(\omega),
			\end{align*}
			Let $t(\omega,E)=min\{t_1,t_2,t_3,2t_4\}$ and
			\begin{align*}
				R_1(\omega)=e^{\nu}[1+\frac{2}{\nu}(r_1(\omega)+r_2(\omega)+R_{f,\nu}(\omega))].
			\end{align*}
			Finally, since 
			\begin{align*}
			S(-t_0,\theta_{t_0} \omega)x_0=x(0,\omega,t_0,x_0)=\{v^\varepsilon(0,\omega,t_0,z(t_0,\omega)u_0),\rho(0,\omega,t_0,\rho_0),\sigma(0,\omega,t_0,\sigma_0)\},
		\end{align*}
		we can obtain that $|S(-t_0,\theta_{t_0} \omega)x_0| \leq 2R_0+R_1(\omega)$ if $
			t_0 \leq t(\omega,E)$, which means $B(0,2R_0+R_1(\omega))$ is a random absorbing set in $H$.
		\end{proof}
	\end{lem}
	Lemma \ref{lem8} gets a random absorbing set in $H$. In order to obtain the asymptotic compactness of the solution in $H$, we need to show that $v^\varepsilon(t)$, $\rho(t)$ and $\sigma(t)$ are bounded in $V$. The following lemma plays a crucial role in the proof of our desired conclusion.
	\begin{lem}\label{lem9}
		There exist a positive number $R_2$ depending on $D$ and universal constants, and a random variable $R_3(\omega)$ depending on the initial data, $\varepsilon_0$, $\nu$, $f$ and universal constants, such that: For every $E > 0$, there exists $t(\omega,E) \leq -1$, such that for all $t_0 \leq t(\omega,E)$ and $x_0=(u_0,\rho_0,\sigma_0)$ $\in{H}^{1}$ with $|x_0| < E$, we have 
		\begin{align}\label{f1}
			\int_{-1}^{0} \|\nabla \sigma(s,\omega,t_0,\sigma_0)\|_{L^2}^{2}ds \leq R_2,
		\end{align}
		\begin{align}\label{f2}
			\int_{-1}^{0} \|\nabla \rho(s,\omega,t_0,\rho_0)\|_{L^2}^{2}ds \leq R_2,
		\end{align}
		\begin{align}\label{f3}
			\int_{-1}^{0} \|\nabla v^\varepsilon(s,\omega,t_0,z(t_0,\omega)u_0)\|_{L^2}^{2}ds \leq R_3(\omega).
		\end{align}
		\begin{proof}
			The differential inequality \eqref{d3} gives 
			\begin{align}\label{d12}
				\frac{d}{dt}( \|\sigma\|_{L^2}^{2}\ + \|\rho\|_{L^2}^{2}\ ) +2D(\|\nabla\rho\|_{L^2}^{2} + \|\nabla\sigma\|_{L^2}^{2}) \leq 0.
			\end{align}
			Integrating \eqref{d12} between $-1$ and $0$ yields
			\begin{align*}
				&\|\sigma(0)\|_{L^2}^{2}-\|\sigma(-1)\|_{L^2}^{2}+\|\rho(0)\|_{L^2}^{2}-\|\rho(-1)\|_{L^2}^{2}+2D\int_{-1}^{0} (\|\nabla \sigma(s)\|_{L^2}^{2}+\|\nabla \rho(s)\|_{L^2}^{2})ds \leq 0.
			\end{align*}
			It follows that
			\begin{align*}
				\int_{-1}^{0} (\|\nabla \sigma(s)\|_{L^2}^{2}+\|\nabla \rho(s)\|_{L^2}^{2})ds &\leq \frac{1}{2D}(\|\sigma(-1)\|_{L^2}^{2}+\|\rho(-1)\|_{L^2}^{2}) \leq \frac{{R_0}^2}{2D}=R_2,
			\end{align*}
			which proves \eqref{f1} and \eqref{f2}.
			Integrating \eqref{d6} between -1 and 0, we get
			\begin{align*}
				\|&v^\varepsilon(0)\|_{L^2}^{2}-\|v^\varepsilon(-1)\|_{L^2}^{2}+\nu \int_{-1}^{0} \|\nabla v^\varepsilon(s)\|_{L^2}^{2}ds\\& \leq \frac{1}{\nu}\int_{-1}^{0} {z^2(s,\omega)}(\|\rho(s)\|_{L^2}^{3}\|\nabla \rho(s)\|_{L^2}+\|f\|_{L^2}^2)ds.
			\end{align*}
			Thus, we have
			\begin{align*}
				\int_{-1}^{0} \|\nabla v^\varepsilon(s)\|_{L^2}^{2}ds \leq \frac{1}{\nu}{R_1}^2(\omega)+\frac{{R_0}^3}{{\nu}^2}(\int_{-1}^{0} z^4(s,\omega) ds +R_2+C_f)=R_3(\omega),
			\end{align*}
			which proves \eqref{f3}.
		\end{proof}
	\end{lem}
	Then, we provide the $H^1$ boundedness of the solution. 
	\begin{lem}\label{lem10}
		There exist two random variables  $R_4(\omega)$ and $R_5(\omega)$ depending on the initial data, $\varepsilon_0$, $\nu$, $D$, $f$ and universal constants, such that for every $E > 0$, there exists $t(\omega,E) \leq -1$, such that for all $t_0 \leq t(\omega,E)$ and $x_0=(u_0,\rho_0,\sigma_0) \in {H}^{1}$  with $|x_0| < E$, we have
		\begin{align}\label{d13}
			\|\nabla v^\varepsilon(0,\omega,t_0,z(t_0,\omega)u_0)\|_{L^2} \leq  R_4(\omega),~for~all~t \in [-1,0],
		\end{align}
		\begin{align}\label{d14}
			\|\nabla \sigma(0,\omega,t_0,\sigma_0)\|_{L^2} \leq R_5(\omega),~for~all~t \in [-1,0],
		\end{align}
		\begin{align}\label{d15}
			\|\nabla \rho(0,\omega,t_0,\rho_0)\|_{L^2} \leq R_5(\omega),~for~all~t \in [-1,0].
		\end{align}
		Therefore, $B(0,R_4(\omega)+2{R_5(\omega)})$ is a random absorbing set in $V$.
		\begin{proof}
			The differential inequality \eqref{b4} gives 
			\begin{align*}
				\frac{d}{dt}\|\nabla v^\varepsilon\|_{L^2}^{2}\ +\nu \|\Delta v^\varepsilon\|_{L^2}^{2} \leq& \frac{C}{{z^2(t,\omega)}}\|\nabla v^\varepsilon\|_{L^2}^{4}+\frac{D}{2}\|\nabla \rho \|_{L^2}^2 +C{z^4(t,\omega)}\|\rho \|_{L^2}^6\\&+C_f+Cz^4(t,\omega),
			\end{align*}
			which implies that 
			\begin{align*}
				\frac{d}{dt}\|\nabla v^\varepsilon\|_{L^2}^{2}\   \leq&  \frac{C}{{z^2(t,\omega)}}\|\nabla v^\varepsilon \|_{L^2}^{2}\|\nabla v^\varepsilon \|_{L^2}^{2}+\frac{D}{2}\|\nabla \rho \|_{L^2}^2+C{z^4(t,\omega)}\|\rho \|_{L^2}^6\\&+C_f+Cz^4(t,\omega).
			\end{align*}
			Let
			\begin{align*}
				f(t)=\frac{C}{{z^2(t,\omega)}}\|\nabla v^\varepsilon(t) \|_{L^2}^{2}.
			\end{align*}
			Applying Lemma \ref{lem1} on an arbitrary interval $[s,0] \subset [-1,0]$, we find
			\begin{align*}
				\|\nabla v^\varepsilon(0)\|_{L^2}^{2} \leq& e^{\int_{s}^{0} f(\xi)d\xi}[\|\nabla v^\varepsilon(s)\|_{L^2}^{2}+\frac{D}{2}\int_{s}^{0}\|\nabla \rho(\xi) \|_{L^2}^{2}d\xi+C\int_{s}^{0} {z^4(\xi,\omega)}\|\rho(\xi) \|_{L^2}^6 d\xi\\
				&+\int_{s}^{0} (C_f+Cz^4(\xi,\omega))d\xi].
			\end{align*}
			Integrating $s$ on the interval $[-1,0]$ yields
			\begin{align*}
				\|\nabla v^\varepsilon(0)\|_{L^2}^{2} \leq& e^{\int_{-1}^{0} f(\xi)d\xi}[\int_{-1}^{0}\|\nabla v^\varepsilon(s)\|_{L^2}^{2}ds+\frac{D}{2}\int_{-1}^{0}\|\nabla \rho(s) \|_{L^2}^2ds+C\int_{-1}^{0} {z^4(\xi,\omega)}\|\rho(\xi) \|_{L^2}^6d\xi\\&+\int_{-1}^{0} (C_f+Cz^4(\xi,\omega))d\xi].
			\end{align*}
			We estimate
			\begin{align}\label{d16}
				\int_{-1}^{0} f(\xi)d\xi &\leq C\mathop {sup}\limits_{t \in [-1,0]} \frac{1}{z^2(t,\omega)} \int_{-1}^{0}\|\nabla v^\varepsilon(\xi) \|_{L^2}^2d\xi \notag\\&\leq CR_3(\omega)\mathop {sup}\limits_{t \in [-1,0]} \frac{1}{z^2(t,\omega)}=r_3(\omega),
			\end{align}
			and
			\begin{align*}
				\int_{-1}^{0} {z^4(\xi,\omega)}\|\rho(\xi) \|_{L^2}^6d\xi &\leq \mathop {sup}\limits_{t \in [-1,0]} z^4(t,\omega)  \mathop {sup}\limits_{t \in [-1,0]}\|\rho(t) \|_{L^2}^6\\& \leq {R_0}^6\mathop {sup}\limits_{t \in [-1,0]} z^4(t,\omega)=r_4(\omega).
			\end{align*}
			We also estimate
			\begin{align*}
				\int_{-1}^{0} (C_f+Cz^4(\xi,\omega))d\xi \leq r_5(f,\omega),
			\end{align*}
			and
			\begin{align*}
				\frac{D}{2}\int_{-1}^{0} \|\nabla \rho(s)\|_{L^2}^{2}ds \leq \frac{D}{2}R_2.
			\end{align*}
			Consequently, let
			\begin{align*}
				R_4(\omega)=e^{r_3(\omega)}[{R_3}^2(\omega)+\frac{D}{2}R_2+r_4(\omega)+r_5(f,\omega)].
			\end{align*}
			This proves \eqref{d13}. Similarly, using \eqref{b26} to \eqref{b32}, we can get the following differential inequality
			\begin{align*}
				\frac{d}{dt}(\|\nabla \rho\|_{L^2}^{2}+\|\nabla \sigma\|_{L^2}^{2}) \leq  \frac{C}{{z^2(t,\omega)}}\|\nabla v^\varepsilon \|_{L^2}^{2} ( \| \nabla \rho\|_{L^2}^{2} + \|\nabla \sigma\|_{L^2}^{2})(\|\nabla \rho\|_{L^2}^{2}+\|\nabla \sigma\|_{L^2}^{2}).
			\end{align*}
			Let
			\begin{align*}
				h(t)=\frac{C}{{z^2(t,\omega)}}\|\nabla v^\varepsilon(t) \|_{L^2}^{2} ( \| \nabla \rho(t)\|_{L^2}^{2} + \|\nabla \sigma(t)\|_{L^2}^{2}).
			\end{align*}
			Applying Lemma \ref{lem1} on an arbitrary interval $[s,0] \subset [-1,0]$ again, we have
			\begin{align*}
				\|\nabla \rho(0)\|_{L^2}^{2}+\|\nabla \sigma(0)\|_{L^2}^{2}  \leq e^{\int_{s}^{0} h(\xi)d\xi}	(\|\nabla \rho(s)\|_{L^2}^{2}+\|\nabla \sigma(s)\|_{L^2}^{2}).
			\end{align*}
			Integrating $s$ on the interval $[-1,0]$ yields
			\begin{align*}
				\|\nabla \rho(0)\|_{L^2}^{2}+\|\nabla \sigma(0)\|_{L^2}^{2}  \leq e^{\int_{-1}^{0} f(\xi)d\xi}\int_{-1}^{0}(\|\nabla \rho(s)\|_{L^2}^{2}ds+\|\nabla \sigma(s)\|_{L^2}^{2})ds.
			\end{align*}
			Lemma \ref{lem6} gives 
			\begin{align*}
				\|\nabla \rho (t)\|_{L^2}^2+\|\nabla \sigma (t)\|_{L^2}^2 \leq  A.
			\end{align*}
			From Lemma \ref{lem9}, we know that
			\begin{align*}
				\int_{-1}^{0}(\|\nabla v^\varepsilon(s)\|_{L^2}^{2}ds \leq R_3(\omega).
			\end{align*}
			Thus we estimate
			\begin{align*}
				\int_{-1}^{0} h(\xi)d\xi \leq CAR_3(\omega)\mathop {sup}\limits_{t \in [-1,0]} \frac{1}{z^2(t,\omega)}=r_6(\omega).
			\end{align*}
			Finally, let
			\begin{align*}
				R_5(\omega)=e^{r_6(\omega)}R_2,
			\end{align*}
			which proves \eqref{d14} and \eqref{d15}. Therefore, we can conclude that $B(0,R_4(\omega)+2{R_5(\omega)})$ is a random absorbing set in $V$.
		\end{proof}
	\end{lem}	
	Now, applying Theorem \ref{thm1}, we are in a position to present the final conclusion in this section, that is, the existence of a compact random attractor in $H$.
	\begin{thm}\label{Theorem7}
		The random dynamical system associated with the stochastic Nernst-Planck-Navier-Stokes equations \eqref{b2} has a compact random attractor in $H$.
	\begin{proof}
		It follows from Lemma \ref{lem10} that there exists a random closed ball in $V$, which absorbs any deterministic bounded sets in $H$. Applying the compactness of embedding $V \hookrightarrow H$, it is clear to deduce that every bounded closed sets in $V$ are asymptotically compact in $H$. Hence, one can easily obtain the existence of a compact absorbing set. With the help of Theorem \ref{thm1}, the desired conclusion is proved. 
		\end{proof}
	\end{thm}
	
	\section{\large Upper Semicontinuity of Random Attractors}\label{sec5}
	In this section, we consider the asymptotic behavior of random attractors for the system \eqref{b2} when $\varepsilon \rightarrow 0$, and establish upper semicontinuity of random attractors for the system \eqref{b2}. To begin with, we introduce the concepts about upper semicontinuity of random attractors in \cite{B.W2014}.
	\begin{thm}\label{Theorem10}
		Let X be a Banach space and $\varphi$ be a dynamical system defined on X. Given $\varepsilon>0$, we suppose that $\varphi_\varepsilon$ is a random dynamical system over a metric system $(\Omega, \mathcal{F}, P, (\theta_t)_{t \in R})$. Suppose the following three statements are satisfied:\par
		(1)for $P$-a.e.$~\omega \in \Omega$, $t \geq 0$, $\varepsilon_n \rightarrow 0$, $x_n$, $x \in X$, $x_n \rightarrow x$, there holds:
		\begin{align}\label{h1}
			\mathop {lim}\limits_{n \rightarrow \infty} \varphi_{\varepsilon_n}(t, \omega)x_n=\varphi(t)x. 
		\end{align}
		\par
		(2)Let $\mathcal{D}$ be a collection of subsets of X. $\mathscr {A}_0$ is the global attractor of $\varphi$ in $X$. Assume that $\varphi_\varepsilon$ has a random attractor $\mathscr{A}_\varepsilon=\{\mathscr{A}_\varepsilon(\omega)\}_{\omega \in \Omega} \in \mathcal{D}$ and a random absorbing set $\mathcal{K}_\varepsilon=\{\mathcal{K}_\varepsilon(\omega)\}_{\omega \in \Omega} \in \mathcal{D}$ such that for some deterministic positive constant c and for $P$-a.e.$~\omega \in \Omega$,
		\begin{align}\label{h2}
			\mathop {lim}\limits_{\varepsilon \rightarrow 0} sup\|\mathcal{K}_\varepsilon(\omega)\|_X \leq c.
		\end{align}	
		\par
		(3)there exists $\varepsilon_0 > 0$ such that for $P$-a.e.$~\omega \in \Omega$, 
		\begin{align}\label{h3}
			\cup_{0<\varepsilon \leq \varepsilon_0} \mathscr{A}_\varepsilon(\omega) ~is ~precompact~ in~ X.
		\end{align}
		Then, $\mathscr{A}_\varepsilon$ is said to be upper semi-continuous, that is, for $P$-a.e.$~\omega \in \Omega$,
		\begin{align*}
			dist_X(\mathscr{A_\varepsilon}(\omega),\mathscr {A}_0) \rightarrow 0,~as~\varepsilon \rightarrow 0.
		\end{align*}
		Here, $dist_X(\cdot,\cdot)$ denotes the Hausdorff semi-distance given by 
		\begin{align*}
			dist_X(A,B) =\mathop {sup}\limits_{a \in A}\mathop {inf}\limits_{b \in B}d(a,b),
		\end{align*}
		for nonempty sets $A,B \subset X$ on a metric space $(X,d)$.
	\end{thm}
	Then, the following Theorem \ref{Theorem11} will be presented to show that condition \eqref{h1} is fulfilled.
	\begin{thm}\label{Theorem11}
		For $0 < \varepsilon \leq 1$, let $u$ and $v^\varepsilon$ be the solutions of systems \eqref{a1} and \eqref{b2}, respectively. Then, for every $\omega \in \Omega$, $t > 0$,
		\begin{align*}
			\mathop {lim}\limits_{\varepsilon \rightarrow 0} \|v^\varepsilon(t,\omega,{v_0}^\varepsilon)-u(t,u_0)\|_{L^2}^2=0.
		\end{align*}
		\begin{proof}
			Let $y^\varepsilon=v^\varepsilon-u$. Then from \eqref{a1} and \eqref{b2}, we find
			\begin{align}\label{e1}
				\frac{dy^\varepsilon(t)}{dt}-\nu \Delta y^\varepsilon(t)=-\frac{1}{z(t,\omega)}B(v^\varepsilon(t))+B(u(t))+(1-z(t,\omega))\rho \nabla \phi+(z(t,\omega)-1)f.
			\end{align}
			Taking the $L^2$ inner product with $y^\varepsilon$ to the equation \eqref{e1}, we get
			\begin{align}\label{add1}
				\frac{1}{2}\frac{d}{dt}\|y^\varepsilon\|_{L^2}^2+\nu \|\nabla y^\varepsilon \|_{L^2}^2=&-e^{\varepsilon \omega(t)}b(v^\varepsilon,v^\varepsilon,y^\varepsilon)+b(u,u,y^\varepsilon)+(1-e^{\varepsilon \omega(t)})\langle \rho \nabla \phi,y^\varepsilon\rangle\notag \\&+(e^{\varepsilon \omega(t)}-1)\langle f,y^\varepsilon \rangle.
			\end{align}
			We can rewrite \eqref{add1} as
			\begin{align*}
				\frac{1}{2}\frac{d}{dt}\|y^\varepsilon\|_{L^2}^2+\nu \|\nabla y^\varepsilon \|_{L^2}^2=&-e^{\varepsilon \omega(t)}b(y^\varepsilon,u,y^\varepsilon)+(1-e^{\varepsilon \omega(t)})b(u,u,y^\varepsilon)+(1-e^{\varepsilon \omega(t)})\langle \rho \nabla \phi,y^\varepsilon\rangle\\&+(e^{\varepsilon \omega(t)}-1)\langle f,y^\varepsilon \rangle.
			\end{align*}
			Applying Holder's, Ladyzhenskaya’s and Young's inequalities, we estimate
			\begin{align*}
				|e^{\varepsilon \omega(t)}b(y^\varepsilon,u,y^\varepsilon)| &\leq Ce^{\varepsilon \omega(t)}\| y^\varepsilon\|_{L^2}\|\nabla y^\varepsilon\|_{L^2}\|\nabla u\|_{L^2}\\ &\leq \frac{\nu}{6}\|\nabla y^\varepsilon \|_{L^2}^2+Ce^{2\varepsilon \omega(t)}\| y^\varepsilon\|_{L^2}^2\|\nabla u\|_{L^2}^2,
			\end{align*}
			and
			\begin{align*}
				|(1-e^{\varepsilon \omega(t)})b(u,u,y^\varepsilon)| &\leq C|1-e^{\varepsilon \omega(t)}|\|u\|_{L^4} \|\nabla u\|_{L^2} \|y^\varepsilon\|_{L^4} \\&\leq C|1-e^{\varepsilon \omega(t)}| \|u\|_{L^2}^{\frac{1}{2}} \|\nabla u\|_{L^2}^{\frac{3}{2}} \|y^\varepsilon\|_{L^2}^{\frac{1}{2}} \|\nabla y^\varepsilon \|_{L^2}^{\frac{1}{2}} \\&\leq C\|\nabla v^\varepsilon\|_{L^2}^2\|y^\varepsilon\|_{L^2}^2+C\|\nabla u\|_{L^2}^2\|y^\varepsilon\|_{L^2}^2+C|1-e^{\varepsilon \omega(t)}|^{\frac{4}{3}}\| u\|_{L^2}^{\frac{2}{3}}\|\nabla u \|_{L^2}^2.
			\end{align*}
			Proceeding as above, we also estimate
			\begin{align*}
				|(1-e^{-\varepsilon \omega(t)})\langle \rho \nabla \Phi,y^\varepsilon \rangle| &\leq |1-e^{-\varepsilon \omega(t)}|\|\rho\|_{L^2}\|\nabla \Phi\|_{L^\infty}\|y^\varepsilon\|_{L^2} \\&\leq |1-e^{-\varepsilon \omega(t)}|\|\rho\|_{L^2}\|\rho \|_{L^4}\|y^\varepsilon\|_{L^2} \\&\leq \frac{3}{2\nu}|1-e^{-\varepsilon \omega(t)}|^2\|\rho\|_{L^3}^2 \| \rho\|_{L^2}^2 +\frac{\nu}{6}\|\nabla y^\varepsilon \|_{L^2}^2,
			\end{align*}
			and
			\begin{align*}
				|(e^{-\varepsilon \omega(t)}-1)\langle  f,y^\varepsilon \rangle| &\leq \frac{3}{2\nu}|1-e^{-\varepsilon \omega(t)}|^2\|f\|_{L^2}^2 +\frac{\nu}{6}\|\nabla y^\varepsilon \|_{L^2}^2.
			\end{align*}	
			Summing the above equalities up, we derive that
			\begin{align*}
				\frac{d}{dt} \|y^\varepsilon\|_{L^2}^2 \leq W_1(t) \|y^\varepsilon\|_{L^2}^2+W_2(t),
			\end{align*}
			where 
			\begin{align*}
				W_1(t)=C[(e^{2\varepsilon \omega(t)}+1) \|\nabla u(t)\|_{L^2}^2+\|\nabla v^\varepsilon(t)\|_{L^2}^2],
			\end{align*}
			and
			\begin{align*}
				W_2(t)=C|1-e^{\varepsilon \omega(t)}|^{\frac{4}{3}}\| u\|_{L^2}^{\frac{2}{3}}\|\nabla u \|_{L^2}^2+\frac{3}{2\nu}|1-e^{-\varepsilon \omega(t)}|^2\|f\|_{L^2}^2+\frac{3}{2\nu}|1-e^{-\varepsilon \omega(t)}|^2\|\rho\|_{L^3}^2 \|\rho\|_{L^2}^2.
			\end{align*}
			Applying Lemma \ref{lem1} on the interval $[0,t]$ yields
			\begin{align*}
				\|v^\varepsilon(t,\omega,{v_0}^\varepsilon)-u(t,u_0)\|_{L^2}^2 \leq (\|e^{-\varepsilon \omega(t)}u_0-u_0\|_{L^2}^2+\int_{0}^{t} W_2(s)ds)e^{\int_{0}^{t} W_1(s)ds}.
			\end{align*}
			We deduce
			\begin{align*}
				\int_{0}^{t} W_1(s)ds \leq C\mathop {sup}\limits_{s \in [0,t]}(e^{2\varepsilon \omega(s)}+1) \int_{0}^{t} \|\nabla u(s)\|_{L^2}^2+\|\nabla v^\varepsilon(s)\|_{L^2}^2 ds  < \infty.
			\end{align*}
			Similarly, we also deduce that
			\begin{align*}
				\int_{0}^{t} W_2(s)ds < \infty.
			\end{align*}	
			Taking limit $\varepsilon \rightarrow 0$, then we obtain
			\begin{align*}
				\mathop {lim}\limits_{\varepsilon \rightarrow 0} \|v^\varepsilon(t,\omega,{v_0}^\varepsilon)-u(t,u_0)\|_{L^2}^2=0,
			\end{align*}
			which completes the proof.		
		\end{proof}
	\end{thm}
	Finally, based on the above theorems, we are ready to establish upper semicontinuity of random attractors with respect to the system \eqref{b2}.
	\begin{thm}\label{Theorem12}
		For $0 < \varepsilon \leq 1$, let $\mathscr{A}(\omega)$ be the global random attractor of system \eqref{a1} and $\mathscr{A_\varepsilon}(\omega)$ be the random attractor of system \eqref{b2}. Then for every $\omega \in \Omega$, $t > 0$,
		\begin{align*}
			\mathop {lim}\limits_{\varepsilon \rightarrow 0} dist_H(\mathscr{A_\varepsilon}(\omega),\mathscr A)=0.
		\end{align*}
		\begin{proof}	 
			Define $\mathcal{K}_\varepsilon(\omega)$ as follows:
			\begin{align*}
				\mathcal{K}_\varepsilon(\omega)=\{u \in H:\|u\|_H^2 \leq \mathcal{M}_\varepsilon(\omega)\},
			\end{align*}
			where $\mathcal{M}_\varepsilon(\omega)$ is given by
			\begin{align*}
				\mathcal{M}_\varepsilon(\omega)=z^{-2}(t,\omega)e^{\nu}[1+\frac{1}{\nu}(r_1(\omega)+r_2(\omega)+R_{f,\nu}(\omega))].
			\end{align*}
			Then for every $0 \leq \varepsilon <1$, $\mathcal{K}_\varepsilon(\omega)$ is a closed absorbing set. For every $\omega \in \Omega$, we have
			\begin{align*}
				\mathop {lim}\limits_{\varepsilon \rightarrow 0}\|\mathcal{K}_\varepsilon(\omega)\|_H^2 \leq \mathop {lim}\limits_{\varepsilon \rightarrow 0}\mathcal{M}_\varepsilon(\omega) =\mathcal{M}_0(\omega).
			\end{align*}
			Then, condition \eqref{h2} is satisfied.
			
			By Lemma \ref{lem10}, $E_\varepsilon(\omega)$ can be defined as: 
			\begin{align*}
				E_\varepsilon(\omega)=\{u \in H:\|u\|_V^2 \leq N_\varepsilon(\omega)\},
			\end{align*}
			where $\mathcal{N}_\varepsilon(\omega)$ is denoted by 
			\begin{align*}
				\mathcal{N}_\varepsilon(\omega)=R_4(\omega)+2R_5(\omega).
			\end{align*}
			Following from the invariance of the random attractor $\mathscr{A_\varepsilon}(\omega)$, we find 
			\begin{align*}
				\cup_{0<\varepsilon \leq 1} \mathscr{A}_\varepsilon(\omega) \subset \cup_{0<\varepsilon \leq 1} E_\varepsilon(\omega),
			\end{align*}
			which implies that 	$\cup_{0<\varepsilon \leq \varepsilon_0} \mathscr{A}_\varepsilon(\omega)$ is bounded in $H^1$. Thanks to the compactness of embedding $H^1 \hookrightarrow L^2$, we conclude that $\cup_{0<\varepsilon \leq \varepsilon_0} \mathscr{A}_\varepsilon(\omega)$ is	precompact in $H$.
			Therefore, by Theorem \ref{Theorem10} along with Theorem \ref{Theorem11}, the proof is completed. 
		\end{proof}
	\end{thm}

	\bibliographystyle{plain}
	
\end{document}